# Estimating the Distribution of Random Parameters in a Diffusion Equation Forward Model for a Transdermal Alcohol Biosensor

Melike Sirlanci[a], Susan E. Luczak[b], Catharine E. Fairbairn[c], Dahyeon Kang[c], Ruoxi Pan[d], Xin Yu[d], and I. G. Rosen[a*]

[a] Modeling and Simulation Laboratory, Department of Mathematics, University of Southern California, E-mail: sirlanci@usc.edu, grosen@math.usc.edu
[b] Department of Psychology, University of Southern California, E-mail: luczak@usc.edu
[c] Department of Psychology, University of Illinois, E-mail: cfairbai@illinois.edu, dkang38@illinois.edu
[d] Department of Mathematics, University of Southern California (undergraduate student)
[*] Corresponding author. Phone: (213)740 2446. Fax: (213)740 2424. E-mail: grosen@math.usc.edu

**Abstract**

We estimate the distribution of random parameters in a distributed parameter model with unbounded input and output for the transdermal transport of ethanol in humans. The model takes the form of a diffusion equation with the input being the blood alcohol concentration and the output being the transdermal alcohol concentration. Our approach is based on the idea of reformulating the underlying dynamical system in such a way that the random parameters are now treated as additional space variables. When the distribution to be estimated is assumed to be defined in terms of a joint density, estimating the distribution is equivalent to estimating the diffusivity in a multi-dimensional diffusion equation and thus well-established finite dimensional approximation schemes, functional analytic based convergence arguments, optimization techniques, and computational methods may all be employed. We use our technique to estimate a bivariate normal distribution based on data for multiple drinking episodes from a single subject.

## 1. Introduction

Researchers and clinicians studying and treating alcohol dependence have long sought the means to continuously and quantitatively monitor blood alcohol levels in naturalistic settings. The ability to do this would be extremely valuable for advancing a wide range of alcohol research and clinical treatment domains, including how alcohol concentrations relate to drinking motives, physical responses to alcohol, behaviors, decision-making, negative consequences, and environmental situational factors over the course of a drinking episode and across drinking episodes. At present the only ways to collect alcohol consumption data in the field are the self-report diary and the breath analyzer, both of which require the active participation of the subject and generally yield inaccurate estimates of blood alcohol level when obtained during naturalistic drinking episodes.

Recently, biosensors that measure transdermal alcohol concentration (TAC), the amount of alcohol diffusing through the skin (Fig. 1.1), have been developed and used, but primarily only as abstinence monitors (for example, in court mandated monitoring of DUI offenders). Because TAC data does not consistently correlate with breath/blood alcohol concentrations (BrAC/BAC) across individuals, devices, and environmental conditions, these devices have not experienced wide-spread acceptance among the research and clinical communities. Indeed, BAC and BrAC are currently, and historically have been, the standard measures of alcohol level intoxication among alcohol researchers and clinicians and unfortunately there is currently no well-established method for producing reliable estimates of BrAC/BAC (eBrAC/eBAC) from TAC data.

The transport and filtering of alcohol by the skin is affected by a number of factors that differ across individuals (e.g., skin layer thickness, porosity, and tortuosity, etc.) and across drinking episodes within individuals (e.g., body and ambient temperature, humidity, subject activity level, skin hydration, vasodilation, etc.). The implication is that, regardless of how reliable and accurate transdermal alcohol device hardware becomes at measuring TAC, the raw TAC data will never consistently map directly onto BrAC/BAC across individuals and drinking episodes.

In our work to date on determining eBrAC/eBAC from TAC (see, for example, [5] and [16]), we have taken a strictly deterministic approach to converting TAC to either BAC or BrAC. First principles physics-based models (a one-dimensional diffusion equation with either infinite or finite speed of propagation and input and output on the boundary) were fit to *individual calibration data* to capture the forward process - the transport of ethanol molecules from the blood, through the skin, and its measurement by the sensor. The result is TAC expressed as a convolution of BAC or BrAC with a kernel or filter. We then deconvolve an estimate of the BAC or BrAC from the biosensor measurements of TAC.

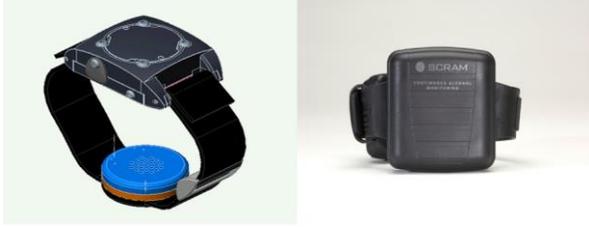

Figure 1.1. (Left Panel) The WRISTAS$^{TM}$ 7 transdermal alcohol biosensor, (Right Panel) The Alcohol Monitoring System (AMS) Secure Continuous Alcohol Monitoring (SCRAM) system.

To greatly reduce the burden on researchers/clinicians and participants/patients and thus significantly increase the feasibility of using these devices, we have been investigating ways to eliminate the need to calibrate the sensor's data analysis system to each individual, each sensor, and across varying current environmental conditions. In particular, we have been investigating the use of our first principles physics/physiological based models to describe the dynamics common to the entire population *(interpreted broadly to include all individuals, devices, and environmental conditions)* and then to attribute all un-modeled sources of uncertainty (primarily due to variations in physiology, hardware, and the environment) observed in individual data to random effects. We refer to this as a *population model*; it takes the form of our earlier deterministic transport models, but with the parameters now being random variables whose distributions are to be estimated based on aggregate population data.

We assume an underlying mathematical framework describing the system dynamics that are common to all individuals, environmental conditions, and devices in the population (e.g., the physics-based model for the transport of ethanol from the blood, through the skin, and measurement by the sensor), but that individual members of the population exhibit variation in the parameters appearing in the model (e.g., the rate at which the alcohol is transported, evaporates, etc.). We then assume that the sensor measures the sum or mean of all of these effects. This is realized in the form of a model based on random partial differential equations and boundary conditions, and then, instead of fitting the unknown parameters in the model directly to individual data, we now estimate the distribution of the random parameters (in the form of probability measures or density functions) based on aggregate population data. In this way, the fit describes the mean behavior of the population.

We consider the (in general, infinite dimensional) discrete time initial value problem given by

$$x_{k+1} = g(t_k, x_k, u_k; q) \quad (1.1)$$
$$x_0 \text{ given}, \quad (1.2)$$

and assume that it describes the dynamics of a process common to the entire population. In addition, it is assumed that we are able to observe some function of the state, or solution, of (1.1), (1.2), $x(u; q) = \{x_k(u, q)\}$, and the input, $u = \{u_k\}$, as given by the output or observation equation

$$y_k(u, q) = h(x_k(u, q), u_k; q). \quad (1.3)$$

In equations (1.1)-(1.3) it is assumed that the parameters $q$ are contained in a set $Q$, the set of admissible parameters, and that their values are specific to a particular individual (interpreted broadly) member of the population. The objective is to estimate the parameters, $q$, or more precisely their distribution, based on population or aggregate data rather than data for a particular individual. We assume $q$ is a random vector with support set $Q$ and that $q \sim \pi(\rho)$, where $\pi(\rho)$ is a family of distributions or push forward measures parameterized by $\rho$ and defined on a sigma field of subsets of $Q$ where $\rho \in \mathcal{R}$, $\mathcal{R}$ the set of feasible parameters. We will assume $\pi(\rho)$ is defined in terms of a joint density function $f(\rho)$; $f(\rho)$ is absolutely continuous, non-negative, and

$$\pi(A; \rho) = P(q \in A) = \int_Q \chi_A d\pi(\rho) = \int_A f(q; \rho) dq,$$

for $A$ an event, where $\rho \in \mathcal{R}$.

The statistical model for the population data upon which the estimation is to be based assumes that the observed data points can be represented by the expected value of the output of the model plus random error. We assume that we have $\nu$ subjects and $\mu_i$ noise-corrupted observations for the i$^{th}$ subject. For $\rho \in \mathcal{R}$, we define

$$v_{i,j}(\rho) = E_{\pi(\rho)}[y_{i,j}(\tilde{u}_i, q)]$$
$$= \int_Q h(x_{i,j}(\tilde{u}_i, q), \tilde{u}_{i,j}; q) d\pi(\rho) \quad (1.4)$$
$$= \int_Q h(x_{i,j}(\tilde{u}_i, q), \tilde{u}_{i,j}; q) f(q; \rho) dq,$$

for $j = 0, 1, \dots, \mu_i$, $i = 1, \dots, \nu$, the mean behavior of the population at time $t_j$ given that the distribution of the parameters $q$ is described by the measure $\pi(\rho)$. In (1.4), $x_i$ and $y_i$ are given by (1.1) – (1.3) with $u = \tilde{u}_i$ and $x_0 = x_{i,0}$. As is typically the case in standard linear regression, we assume our observations are given by

$$V_{i,j} = v_{i,j}(\rho_0) + \varepsilon_{i,j} = E_{\pi(\rho_0)}[y_{i,j}(\tilde{u}_i, q)] + \varepsilon_{i,j}$$
$$= \tilde{y}_{i,j}(\tilde{u}_i, \rho_0), \; j = 0, 1, \dots, \mu_i, \; i = 1, \dots, \nu, \quad (1.5)$$

where the $\varepsilon_{i,j}$, $j = 1, \dots, \mu_i$, $i = 1, \dots, \nu$, represent measurement noise which are assumed to be independent, identically-distributed random variables with mean 0 and common variance $\sigma^2$, and $\rho_0 \in \mathcal{R}$ represents the *true* values of the parameters.

Our estimation problem minimizes prediction error based on a least squares (or MLE if $\varepsilon_{i,j} \sim N(0, \sigma^2)$) approach. For $V = \{V_{i,j}\}$, our estimator $\hat{\rho}$ is given by

$$\hat{\rho} = arg \min_{\rho \in \mathcal{R}} J(\rho; V)$$
$$= arg \min_{\rho \in \mathcal{R}} \{\sum_{i=1}^{\nu} J_{i,\mu}(\rho; V) + \lambda(\rho)\} \quad (1.6)$$
$$= arg \min_{\rho \in \mathcal{R}} \{\sum_{i=1}^{\nu} \sum_{j=0}^{\mu_i} \alpha_{i,j} |v_{i,j}(\rho) - V_{i,j}|^2 + \lambda(\rho)\}$$

where $\lambda(\rho)$ is a regularization term, $\mu = \{\mu_i\}_{i=1}^{\nu}$, the $\alpha_{i,j}$ are positive weights, $v_{i,j}(\rho)$ is given by (1.4) and $V = \{V_{i,j}\}$ by (1.5). In general (and in fact, in the problem of interest to us here), the fact that the dynamical system (1.1)-(1.3) may be (or is) infinite dimensional, means that the optimization problem given in (1.6) is not directly amenable to computation. Consequently, when we apply these ideas to the alcohol biosensor problem in Section 2 below, we replace the infinite dimensional state equation with a finite dimensional approximation. Relevant convergence issues which are the central focus of our research, are discussed in Section 4. We note that if one is trying to estimate the shape of $f$ and the parameters $\rho$ are the coefficients of basis elements with the dimension of the feasible parameter space $\mathcal{R}$ increasing in dimension in order to achieve some level of convergence, then the presence of regularization in the form of the term $\lambda(\rho)$, would likely be essential. In our studies here, we are not estimating the shape of $f$ and the dimension of $\mathcal{R}$ is fixed and relatively low (~10), and hence we found that regularization was unnecessary. Consequently we do not include $\lambda(\rho)$ (i.e. we assume $\lambda(\rho) = 0$) in the remainder of our treatment below. It would not however, be difficult to include it, and we may do that in future work where the focus is either the estimation of the shape of the density or the corresponding measure, $\pi$, as in [4].

Our problem as posed is a classic mixed effects model (see, for example, [6,7], and the numerous references contained therein) featuring both inter- and intra-subject variation and uncertainty. The standard statistical model requires a function, which describes the relationship between the observations and the design variables of the form

$y_{ij} = f(x_{ij}, q_i) + \varepsilon_{ij}, \quad i = 1, \ldots, \nu, \quad j = 1, \ldots, \mu_i,$

where $y_{ij}$ is the j$^{\text{th}}$ observation from the i$^{\text{th}}$ individual, the $x'_{ij}s$ are the known design or state variables, the $q_i$'s are the parameters specific to i$^{\text{th}}$ individual, and $\varepsilon_{ij}$ (typically $\sim N(0, \sigma^2)$) are i.i.d. with mean zero and common variance $\sigma^2$. Moreover, it is also assumed that $q_i = h(c_i, \mu, \eta_i)$ where $c_i's$ are the known vectors of the covariates, $\mu$ is the unknown vector of fixed effects, and the $\eta_i$ are random, for example, $\eta_i \sim N(0, \Omega)$. In our case here, in this model, the parameters to be estimated are $\mu, \sigma, \Omega$ with $\Omega$ the covariance matrix. When f and $h$ are linear, estimation of these parameters is relatively straight forward using standard approaches. However, when f is nonlinear, or even more complex, for example, as in our case here wherein, due to the PDE-based model, a closed form finite dimensional representation is unavailable, and there is a clear lack of independence in the observations, a number of significant statistical and computational challenges result. Indeed, the state variables, $x_{ij}$, in our model are the solution of a partial differential equation and consequently the representation for the function $f$ is in terms of a semigroup of infinite dimensional operators on a Hilbert space whose dependence on the parameters, $q_i$, is highly nonlinear. As a result the use of any estimation technique that relies on a likelihood (e.g. MLE, EM, Stochastic Approximation Expectation Minimization, Bayesian estimation with MCMC, etc.) can be daunting (see, for example, [10]). Indeed, such an approach would in one way or another, likely involve repeated simulation requiring repeated solution of the PDE which is something we are trying to avoid. Another method involves estimating the state from the observations, then using nonlinear regression together with the PDE to estimate the parameters. In our case we have a one dimensional measurement of an infinite dimensional state and observability could very likely be an issue. Also, if one opted to use the Kalman filter to estimate the state based on the observations, its implementation is essentially equivalent to solving the PDE and would itself require finite dimensional approximation. In general, this approach often yields inaccurate estimates.

Along with measurements across different subjects, we have longitudinal measurements for each subject, which of course one would expect to be dependent. In addition, from a pharmacokinetic (PK) standpoint, our definition of "subject" or population is somewhat non-standard in that it refers to not only individual participants and their various physiological differences, but also environmental conditions such as temperature and humidity, etc. as well as hardware related uncertainties. Thus, in our model, the $q_i's$ describe unmodeled phenomena present both across different individuals and within the same individual. Indeed our working hypothesis is that each observation represents a sum or average of any number of diffusion processes all at work simultaneously. In addition, because human subjects are involved, it would be both costly and time consuming to collect simultaneous BrAC and TAC data from enough individuals wearing enough sensors and under sufficiently varied environmental conditions to estimate the distribution of the $q_i's$ directly. Consequently, although it can have bias problems [14], it seemed most appropriate that our estimator in (1.6) take the form of what is known in the PK literature as the *naïve pooled data estimator*.

Our general approach relies to some extent on two relatively recent papers: 1) Banks and Thompson's framework for the estimation of probability measures in random abstract evolution equations and the convergence of finite dimensional approximations in the Prohorov metric [4], and 2) Gittelson, Andreev, and Schwab's theory for random abstract parabolic partial differential equations with dynamics defined in terms of coercive sesquilinear forms [9]. The approach in [9] is novel in the way that it treats the random variable as another "space-like" independent variable in the PDE. In this way, finite dimensional approximation is handled in much the same way that it is for the standard deterministic space variables. We use the framework in [9] together with the generation and approximation results from linear semigroup theory, (i.e. the Hille-Yosida-Phillips and Trotter Kato theorems, see, for example, [3,11,15]) to establish that the sufficient conditions for a Banks–Thompson-like convergence result are satisfied in the case of regularly dissipative systems with random parameters whose distributions can be described by appropriately parameterized probability density functions.

For a number of reasons, the approach we take here, and in particular our population model as defined in Section 3 below, is especially well suited for our estimation problem as given in (1.6) with the underlying dynamics (1.1), (1.2) being described by a random PDE. Indeed, 1) it does not require repeated simulation, 2) it takes particular advantage of the underlying parabolic structure of the model's state equation, 3) it lends itself extremely well to functional analytic arguments for convergence of the estimators based on finite dimensional approximation, the central focus of this study, 4) based on our working hypothesis concerning our data as stated in the previous paragraph and the statistical model given in (1.5), the output of the population model is precisely what is required to evaluate the naïve pooled data based performance index, $J(\rho; V)$ given in (1.6), and 5) it is especially well suited for deconvolving an estimate of BrAC from TAC which is our ultimate goal [20].

In this paper, we focus solely on the application of our framework to the alcohol biosensor problem outlined above.

The abstract functional analytic approximation and convergence theory on which the results presented here are based was established in [18] and [19]. In addition, in our treatment here, we are only concerned with the fitting of the random parameters in the forward model. The inverse problem involving the deconvolution of the BAC/BrAC from the TAC signal once the forward model with random elements has been fit is treated elsewhere [20].

An outline of the remainder of the paper is as follows. In Section 2 we discuss our diffusion based distributed parameter model for the transdermal transport of ethanol, its abstract formulation as an infinite dimensional dynamical system with unbounded input and output, and the corresponding discrete-time input-output system on which our general framework is based. In Section 3 we discuss the treatment of initial value problems involving regularly dissipative operators with random coefficients as in [9] and [17], and apply it to the system of interest to us here discussed in Section 2. Section 4 outlines our finite dimensional approximation and convergence results, in Section 5 we discuss a consistency result for our estimator, and in Section 6 we present and discuss our numerical results using actual experimental/clinical data for the alcohol biosensor problem. Section 7 contains a brief discussion of future research, and some concluding remarks.

We use standard notation throughout. For example, we denote the space of square Lebesgue or Bochner integrable functions defined on an interval $(a,b)$ with range in the normed linear space $X$ by $L_2(a,b;X)$, and we use $C(a,b;X)$ when the functions are continuous. When $X = \mathbb{R}$, the range space is omitted. For normed linear spaces $X$ and $Y$, $\mathcal{L}(X,Y)$ denotes the space of bounded (continuous) linear operators defined on $X$ with range in $Y$. Unless explicitly stated otherwise, all Hilbert space norms are the ones induced by the standard inner product on that space. We occasionally use "dot" notation to denote weak or strong derivatives with respect to time.

## 2. A Distributed Parameter Model for a Transdermal Alcohol Biosensor and its Abstract Formulation

Our ethanol transport model is based on diffusion, or Fick's law, and consequently it is described by abstract parabolic operators. When formulated abstractly in a Gelfand triple setting, these operators are examples of what are known as regularly dissipative operators and can be shown to generate analytic or holomorphic semigroups (see [11,15,22]).

After converting to what are essentially dimensionless quantities (see [16]), we obtain the input/output model

$$\frac{\partial \varphi}{\partial t}(t,\eta) = q_1 \frac{\partial^2 \varphi}{\partial \eta^2}(t,\eta), 0 < \eta < 1, \ t > 0 \quad (2.1)$$

$$q_1 \frac{\partial \varphi}{\partial \eta}(t,0) - \varphi(t,0) = 0, \ t > 0, \quad (2.2)$$

$$q_1 \frac{\partial \varphi}{\partial \eta}(t,1) = q_2 u(t), \ t > 0 \quad (2.3)$$

$$\varphi(0,\eta) = \varphi_0, \ 0 < \eta < 1 \quad (2.4)$$

$$y(t) = \varphi(t,0), \ t > 0, \quad (2.5)$$

in the form of an initial-boundary value problem for a one dimensional diffusion equation with input and output on the boundary and two unknown parameters, $q = (q_1, q_2)$. In the system (2.1) – (2.5) $\varphi(t,\eta)$ is essentially the concentration of ethanol in the interstitial fluid in the epidermal layer of the skin at depth $\eta$ and time t, $u$ is the concentration of alcohol in the blood (BAC) as measured by a breath analyzer (BrAC), and y is the (TAC). The boundary condition (2.2) models the evaporation of ethanol at the skin surface, condition (2.3) captures the exchange of ethanol molecules between the (blood fed) dermal and epidermal layers of the skin. The output equation (2.5) models the biosensor measured TAC at the skin surface. We assume that there is no alcohol in the skin initially, so in general $\varphi_0 = 0$ in (2.4).

Using the tools of functional analysis and linear semigroup theory, we reformulate (2.1) – (2.5) as a discrete-time SISO system with state space an infinite dimensional Hilbert space. In (2.1)-(2.5) the input and output are on the boundary and consequently the resulting continuous time input and output operators are unbounded with respect to the usual state space for such a system, $L_2(0,1)$. However, in the discrete or sampled time formulation which is of primary interest to us here, they become bounded.

Let V and H be Hilbert spaces with the embeddings V $\hookrightarrow$ H $\hookrightarrow$ V* dense and continuous, where V* denotes the space of continuous linear functionals on V. Let $\langle \cdot,\cdot \rangle$ and $|\cdot|$ denote the H inner product and norm, respectively, and let $\|\cdot\|$, denote the norm on V. For q $\in \{Q, d_Q\}$, a compact metric space, let $a(q;\cdot,\cdot) : V \times V \to \mathbb{R}$ be a bilinear form satisfying the following three conditions:

i. (Boundedness) $|a(q;\psi_1,\psi_2)| \leq \alpha_0 \|\psi_1\|\|\psi_2\|$,
    $\psi_1, \psi_2 \in V, q \in Q$,

ii. (Coercivity) $a(q;\psi,\psi) + \lambda_0|\psi|^2 \geq \mu_0\|\psi\|^2$,
    $\psi \in V, q \in Q$

iii. (Measurability) For $\psi_1, \psi_2 \in V$, the map $q \to a(q;\psi_1,\psi_2)$ is measurable with respect to all measures $\pi(\rho), \rho \in \mathcal{R}$.

For q $\in \{Q, d_Q\}$ let b(q),c(q) $\in V^*$, and consider an input/output system in weak form as given by

$$\langle \dot{\varphi},\psi \rangle_{V^*,V} + a(q; \varphi,\psi) = \langle b(q),\psi \rangle_{V^*,V} u,$$
$$y = \langle c(q),\varphi \rangle_{V^*,V}, \quad \psi, \varphi \in V, \quad (2.6)$$

where $\varphi(0) = \varphi_0, \in H$ and $\langle \cdot,\cdot \rangle_{V^*,V}$ denotes the natural extension of the H inner product to the duality pairing between V and V*. If we set $W(0,T) = \{\psi: \psi \in L_2(0,T;V), \dot{\psi} \in L_2(0,T;V^*)\}$ and u $\in L_2(0,T)$ it can be shown (see, for example, [13]) that the system (2.6) admits a unique solution $\varphi \in W(0,T)$ that depends continuously on u $\in L_2(0,T)$. It follows that $W(0,T) \subseteq C(0,T,H)$ and that y $\in L_2(0,T)$.

For q $\in Q$, the q-dependent bilinear form on V $\times$ V, $a(q;\cdot,\cdot): V \times V \to \mathbb{R}$, defines a bounded linear operator $A(q) \in \mathcal{L}(V,V^*)$ by $\langle A(q)\psi_1,\psi_2 \rangle_{V^*,V} = -a(q;\psi_1,\psi_2)$, for $\psi_1, \psi_2 \in V$. Then, if we let $\mathcal{H}$ denote any of the spaces V, H or V*, we can consider the linear operator A(q) to be the unbounded linear operator, $A(q): D_q \subset \mathcal{H} \to \mathcal{H}$ where $D_q = V$ in the case $\mathcal{H} = V^*$, and $D_q = \{\psi \in V: A(q)\psi \in \mathcal{H}\}$ in the case $\mathcal{H} = H$ or $\mathcal{H} = V$. It can then be shown (see, for example, [2,3,22]) that $A(q)$ is a closed, densely defined unbounded linear operator on $\mathcal{H}$ and it is the infinitesimal generator of an analytic semigroup of bounded linear operators, $\{e^{A(q)t}: t \geq 0\}$ on $\mathcal{H}$.

For $q \in Q$, define the bounded linear operators $B(q): \mathbb{R} \to V^*$ and $C(q): V \to \mathbb{R}$ by $\langle B(q)u, \psi \rangle_{V^*,V} = \langle b(q), \psi \rangle_{V^*,V} u$, and $C(q)\psi = \langle c(q), \psi \rangle_{V^*,V}$, respectively, for $\psi \in V$ and $u \in \mathbb{R}$. The input/output system can now be written *formally* in the standard state space form as

$$\dot{x}(t) = A(q)x(t) + B(q)u(t),$$
$$y(t) = C(q)x(t), \; t > 0, \; x(0) = x_0 \in H \quad (2.7)$$

where the state $x(t) = \varphi(t, \cdot)$. Using the fact that $\{e^{A(q)t}: t \geq 0\}$ is an analytic semigroup on the spaces $V, H$ and $V^*$, it follows that a so called mild solution ([11,15]) to the state equation in (2.7) x is given by

$$x(t; q) = e^{A(q)t}x_0 + \int_0^t e^{A(q)(t-s)}B(q)u(s)\,ds, t \geq 0, (2.8)$$

where $x$ in (2.8) is in $W(0, T)$.

Now let a sampling time $\tau > 0$ and $x_0 \in V$ be given, and consider zero order hold inputs of the form $u(t) = u_i$, $t \in [i\tau, (i+1)\tau)$, $i = 0,1,2,\ldots$. Then, under the assumptions we have made here thus far (see [19]), it can be shown that using (2.8), it follows from the properties of analytic semigroups generated by regularly dissipative operators that $x_{i+1} = \hat{A}(q)x_i + \hat{B}(q)u_i$, $y_i = \hat{C}(q)x_i$, $i = 0, 1, 2, \ldots$, (2.9) where $\hat{A}(q) = e^{A(q)\tau} \in \mathcal{L}(V,V)$, $\hat{B}(q) = \int_0^\tau e^{A(q)s}B(q)ds \in \mathcal{L}(\mathbb{R},V)$, $\hat{C}(q) = C(q) = c(q) \in \mathcal{L}(V,\mathbb{R}) = V^*$, $x_i = x(i\tau)$ and $y_i = y(i\tau)$, $i = 0, 1, 2, \ldots$.

Boundedness of the operators $\hat{A}(q)$ and $\hat{B}(q)$ in (2.9) follows once again from the fact that $\{e^{A(q)t}: t \geq 0\}$ is an analytic semigroup on $V, H$ and $V^*$ ([2,3,13,22]). Indeed, the coercivity assumption, Assumption (ii) (possibly together with a change of variables), implies that we may assume without loss of generality that the operator $A(q)$ from either $H$ into $H$ or $V$ into $V^*$, is invertible with bounded inverse. Consequently, it follows that

$$\hat{B}(q) = \int_0^\tau e^{A(q)s}B(q)ds = A(q)^{-1}e^{A(q)s}B(q)\Big|_0^\tau$$
$$= (\hat{A}(q) - I)A(q)^{-1}B(q)$$
$$= (\hat{A}(q) - I)A(q)^{-1}b(q) \in \mathcal{L}(V,\mathbb{R}).$$

Let $Q$ be a closed and bounded subset of $\mathbb{R}^2$ endowed with the Euclidean metric, let $H = L_2(0,1)$ together with the standard inner product $\langle \psi_1, \psi_2 \rangle = \int_0^1 \psi_1(x)\psi_2(x)\,dx$, and norm denoted by $|\cdot|$, and let $V$ be the Sobolev space $V = H_1(0,1)$ together with its standard inner product $\langle\langle \psi_1, \psi_2 \rangle\rangle = \int_0^1 \psi_1(x)\psi_2(x)\,dx + \int_0^1 \psi_1'(x)\psi_2'(x)\,dx$ and norm denoted by $\|\cdot\|$. Then we have the usual dense and continuous embeddings $V \hookrightarrow H \hookrightarrow V^*$, where $V^*$ denotes the space of distributions dual to $V$. The forms and functions $a(q;\cdot,\cdot): V \times V \to \mathbb{R}$, $b(q;\cdot): V \to \mathbb{R}$ and $c(\cdot): V \to \mathbb{R}$ are given by

$$a(q; \psi_1, \psi_2) = \psi_1(0)\psi_2(0) \quad (2.10)$$
$$+ q_1 \int_0^1 \psi_1'(x)\psi_2'(x)\,dx, \psi_1, \psi_2 \in V$$

$\langle b(q), \psi \rangle_{V^*,V} = q_2 \psi(1)$, and $\langle c(q), \psi \rangle_{V^*,V} = \psi(0)$, for $\psi \in V$. It follows that $b(q) = q_2\delta(\cdot - 1) \in V^*$ and $c(q) = \delta \in V^*$, where $\delta$ denotes the Dirac delta distribution, or unit impulse at zero. It is not difficult to argue that Assumptions (i)-(iii) hold for the form $a(\cdot;\cdot,\cdot)$ as given in (2.10) above. See [19] for a more abstract, detailed and rigorous description of how we deal with input signals on the boundary of the domain.

## 3. Systems Governed by Regularly Dissipative Operators with Random Parameters

In this section, we use ideas from [9] to consider systems of the form (2.1)-(2.5) with the parameters $q \in Q$ random. Let $q$ be a $p$ dimensional random vector with support in $\prod_{i=1}^p [a_i, b_i]$ where $-\infty < \bar{\alpha} < a_i < b_i < \bar{\beta} < \infty$. Let $\vec{a}, \vec{b} \in \mathbb{R}^p$ be given by $\vec{a}=[a_i]$ and $\vec{b}=[b_i]$, and let $\Theta$ be a parameter set that is a compact subset of $\mathbb{R}^r$ for some $r$. We assume that $q$ has distribution described by the absolutely continuous cdf, $F(q; \vec{a}, \vec{b}, \theta)$, or equivalently by the push forward measure $\pi = \pi(\vec{a}, \vec{b}, \theta)$, where $\theta \in \Theta$.

Let $a(q;\cdot,\cdot)$ denote a sesquilinear form on $V \times V$ satisfying the conditions (i) - (iii) given in Section 2, and in particular that the function $q \mapsto a(q; v, w)$ is $\pi-$measurable for any $v, w \in V$. By appropriately defining new function spaces, it is possible to embed the randomness in the sesquilinear form $a(q;\cdot,\cdot)$, or equivalently in the operator, $A(q)$, in (2.7) into these spaces. Consequently, the input output system (2.7) can be stated in a way that makes the stochasticity in the operators, the state and the output effectively invisible and thus amenable to analysis and approximation using standard (deterministic) linear semigroup theory. In effect, the random variables are treated the same way as the space variables in the underlying PDE.

Toward this end, we define $\mathcal{V} = L_\pi^2(Q;V)$ and $\mathcal{H} = L_\pi^2(Q;H)$. It then follows that $\mathcal{V}, \mathcal{H}, \mathcal{V}^*$ form a Gelfand triple of separable Hilbert spaces with $\mathcal{V} \subseteq \mathcal{H} \subseteq \mathcal{V}^*$ by identifying $\mathcal{H}$ with its dual $\mathcal{H}^*$ and identifying $\mathcal{V}^*$ with the Bochner space $L_\pi^2(Q;V^*)$. We also define the $\pi$-averaged sesquilinear form $a(\cdot,\cdot): \mathcal{V} \times \mathcal{V} \to \mathbb{C}$ by
$$a(v,w) = \int_Q a(q; v(q), w(q))d\pi =$$
$$E_\pi[a(q; v(q), w(q))] \quad (3.1)$$
where, $v, w \in \mathcal{V}$. Assumptions (i) and (ii) guarantee that $a(\cdot,\cdot)$ given in (3.1) is a bounded and coercive sesquilinear form on $\times \mathcal{V}$, and therefore that $-a$ induces a bounded linear map $\mathcal{A}$ from $\mathcal{V}$ into $\mathcal{V}^*$. It follows (see, once again, for example, [2,3,22]) that the operator $\mathcal{A}: \text{Dom}(\mathcal{A}) \subseteq \mathcal{H} \to \mathcal{H}$ where $\text{Dom}(\mathcal{A}) = \{v \in \mathcal{V}: \mathcal{A}v \in \mathcal{H}\}$ is the infinitesimal generator of an analytic semigroup of bounded linear operators $\mathcal{T} = \{\mathcal{T}(t): t \geq 0\}$ on $\mathcal{H}$, and, moreover, that $\mathcal{T}$ can be extended to an analytic semigroup on $\mathcal{V}^*$ and restricted to an analytic semigroup on $\mathcal{V}$.

We assume that the maps $q \mapsto \langle b(q), \psi(q)\rangle_{V^*,V}$, and $q \mapsto \langle c(q), \psi(q)\rangle_{V^*,V}$ are $\pi$-measurable for any $\psi \in \mathcal{V}$, and that $\|b(q)\|_{V^*}$ and $\|c(q)\|_{V^*}$ are uniformly bounded for $q \in Q$. (Assuming that $b, c \in \mathcal{V}^*$ would be fine as well.) We then define the two bounded linear operators $\mathcal{B}: \mathbb{R} \to \mathcal{V}^*$ and $\mathcal{C}: \mathcal{V} \to \mathbb{R}$ by the expressions $\langle \mathcal{B}u, v\rangle = \int_Q \langle b(q), v(q)\rangle_{V^*,V} d\pi(q) u = E_\pi[\langle b(q), v(q)\rangle_{V^*,V}]u$ and $\mathcal{C}v = \int_Q \langle c(q), v(q)\rangle_{V^*,V} d\pi(q) = E_\pi[\langle c(q), v(q)\rangle_{V^*,V}]$, respectively, for $v \in \mathcal{V}$ and $u \in \mathbb{R}$.

We then consider the continuous time input/output system set in the spaces $\mathcal{V}, \mathcal{H}, \mathcal{V}^*$, given by
$$\dot{x}(t) = \mathcal{A}x(t) + \mathcal{B}u(t), \quad (3.2)$$
$$\mathcal{Y}(t) = \mathcal{C}x(t), \ t > 0, \ x(0) = x_0 \in \mathcal{H}. \quad (3.3)$$

The mild solution to the initial value problem given in (3.2), (3.3) is then given by
$$x(t) = \mathcal{T}(t)x_0 + \int_0^t \mathcal{T}(t-s)\mathcal{B}u(s)\,ds, \ t \geq 0, \quad (3.4)$$
and therefore, from (3.3), for $t \geq 0$, that
$$\mathcal{Y}(t) = \mathcal{C}\mathcal{T}(t)x_0 + \int_0^t \mathcal{C}\mathcal{T}(t-s)\mathcal{B}u(s)\,ds, \quad (3.5)$$

Once again let the sampling time $\tau > 0$ be given and consider zero order hold inputs of the form $u(t) = u_i$, $t \in [i\tau, (i+1)\tau)$, $i = 0, 1, 2, \ldots$. Set $x_i = x(i\tau)$ and let $\mathcal{Y}_i = \mathcal{Y}(i\tau)$, $i = 0, 1, 2, \ldots$. It then follows from the variation of parameters formula for systems governed by analytic semigroups, (3.4), and (3.5) that
$$x_{j+1} = \hat{\mathcal{A}}x_j + \hat{\mathcal{B}}u_j, \ \mathcal{Y}_j = \hat{\mathcal{C}}j, \ j = 0, 1, 2, \ldots \quad (3.6)$$
with $x_0 \in \mathcal{V}$, where $\hat{\mathcal{A}} = \mathcal{T}(\tau) \in \mathcal{L}(\mathcal{V}, \mathcal{V})$, $\hat{\mathcal{B}} = \int_0^\tau \mathcal{T}(s)\mathcal{B}ds \in \mathcal{L}(\mathbb{R}, \mathcal{V})$, and $\hat{\mathcal{C}} = \mathcal{C} \in \mathcal{L}(\mathcal{V}, \mathbb{R}) = \mathcal{V}^*$. Boundedness of the operators $\hat{\mathcal{A}}$ and $\hat{\mathcal{B}}$ follows from the fact that $\{\mathcal{T}(t): t \geq 0\}$ is an analytic semigroup on $\mathcal{V}, \mathcal{H}$ and $\mathcal{V}^*$ ([2,3,13,22]). Once again, without loss of generality, we may assume that $\mathcal{A}: \mathrm{Dom}(\mathcal{A}) \subset \mathcal{V}^* \to \mathcal{V}^*$ is invertible with bounded inverse, from which it follows that $\hat{\mathcal{B}} = \int_0^\tau \mathcal{T}(s)\mathcal{B}ds = \mathcal{A}^{-1}\mathcal{T}(s)\mathcal{B}|_0^\tau = (\hat{\mathcal{A}} - I)\mathcal{A}^{-1}\mathcal{B} \in \mathcal{L}(\mathbb{R}, \mathcal{V}) = \mathcal{V}^*$.

It can be shown ([9,17]) that $x(t)$ given by (3.4) agrees $\pi$-almost everywhere with $x(t)$ given in (2.8) for all $t \geq 0$, and consequently it follows that $\mathcal{Y}(t) = \mathcal{C}x(t) = E_\pi[y(t; q)] = E_\pi[C(q)x(t; q)]$, for all $t \geq 0$, and therefore that
$$\mathcal{Y}_i = \hat{\mathcal{C}}x_i = E_\pi[y_i(q)] = E_\pi[\hat{C}(q)x_i(q)]. \quad (3.7)$$
We refer to (3.2),(3.3) or (3.6),(3.7) as our *population model*.

## 4. Finite Dimensional Approximation, Convergence and Computational Considerations

In light of the final expression in Section 3, we may formulate our estimation/optimization problem (1.6) as follows:

($\mathcal{P}$) Given $\nu$ data sets $(\tilde{u}_i, \tilde{y}_i) = \left(\{\tilde{u}_{i,j}\}_{j=0}^{\mu_i - 1}, \{\tilde{y}_{i,j}\}_{j=0}^{\mu_i}\right)_{i=1}^\nu$, determine $\vec{a}^* = [a_i^*]$, $\vec{b}^* = [b_i^*]$ and $\theta^* \in \Theta$, feasible, which minimize
$$J(\vec{a}, \vec{b}, \theta) = \sum_{i=1}^\nu \sum_{j=0}^{\mu_i} |\mathcal{Y}_j(\tilde{u}_i, \vec{a}, \vec{b}, \theta) - \tilde{y}_{i,j}|^2,$$
where $\{\mathcal{Y}_j(\tilde{u}_i, \vec{a}, \vec{b}, \theta)\}_{j=0}^{\mu_i}$ is given as in (3.6) with $u_j = \tilde{u}_{i,j}$, $j = 0, 1, 2, \ldots, \mu_i - 1$, and $i = 1, 2, \ldots, \nu$. ∎

Henceforth we assume that the measures $\pi = \pi(\vec{a}, \vec{b}, \theta)$ are described by a family of joint density functions $f = f(\vec{a}, \vec{b}, \theta)$, and let $\bar{Q} = \prod_{i=1}^p [\bar{\alpha}, \bar{\beta}]$. We will require the following assumptions on the family of densities $f(\vec{a}, \vec{b}, \theta)$:

iv. The maps $(\vec{a}, \vec{b}, \theta) \mapsto f(q; \vec{a}, \vec{b}, \theta)$ are continuous on $\mathbb{R}^p \times \mathbb{R}^p \times \Theta$ for $\pi$-almost every $q \in \bar{Q} = \prod_{i=1}^p [\bar{\alpha}, \bar{\beta}]$.

v. There exist positive constants $\gamma, \delta$ such that
$$0 < \gamma \leq f(q; \vec{a}, \vec{b}, \theta) \leq \delta < \infty,$$
for $\pi$-almost every $q \in \bar{Q} = \prod_{i=1}^p [\bar{\alpha}, \bar{\beta}]$.

Assumptions (i)-(v) are sufficient to establish that the maps $(\vec{a}, \vec{b}, \theta) \mapsto \mathcal{Y}_k(\tilde{u}_i, \vec{a}, \vec{b}, \theta)$ are continuous for $k = 0, 1, 2, \ldots, \mu_i - 1$, and $i = 1, 2, \ldots, \nu$, and therefore that the map $(\vec{a}, \vec{b}, \theta) \mapsto J(\vec{a}, \vec{b}, \theta)$ is continuous. Consequently it follows from compactness that problem ($\mathcal{P}$) has a solution.

Solving problem ($\mathcal{P}$) requires finite dimensional approximation. For each $N = 1, 2, \ldots$, let $\bar{\mathcal{H}} = L_\pi^2(\bar{Q}; H)$, and $\bar{\mathcal{V}} = L_\pi^2(\bar{Q}; V)$. Let $\vec{a}^N = [a_i^N]$ and $\vec{b}^N = [b_i^N]$ be vectors in $\mathbb{R}^p$ with $-\infty < \bar{\alpha} < a_i^N < b_i^N < \bar{\beta} < \infty$, and let $\theta^N \in \Theta$, and set $Q^N = \prod_{i=1}^p [a_i^N, b_i^N]$, $\mathcal{H}^N = L_\pi^2(Q^N; H)$, and $\mathcal{V}^N = L_\pi^2(Q^N; V)$. Let $\mathcal{U}^N$ be a finite dimensional subspace of $\mathcal{V}^N$. Let $\mathcal{J}^N: \bar{\mathcal{H}} \to \mathcal{H}^N$ be such that $\mathrm{Im}(\mathcal{J}^N) = \mathcal{H}^N$ and $|\mathcal{J}^N x|_{\mathcal{H}^N} \leq |x|_{\bar{\mathcal{H}}}$. Let $\mathcal{P}^N: \mathcal{H}^N \to \mathcal{U}^N$ be the orthogonal projection of $\mathcal{H}^N$ onto $\mathcal{U}^N$ and define $\mathcal{J}^N := \mathcal{P}^N \circ \mathcal{J}^N$. Define $\mathcal{A}^N: \mathcal{U}^N \to \mathcal{U}^N$ by
$$\langle \mathcal{A}^N v^N, w^N\rangle = -a(v^N, w^N)$$
$$= -\int_{Q^N} a(q; v^N(q), w^N(q)) d\pi(q)$$
$$= -\int_{Q^N} a(q; v^N(q), w^N(q)) f(q; \vec{a}^N, \vec{b}^N, \theta^N) dq,$$
where, $v^N, w^N \in \mathcal{U}^N$.

Using assumptions (i)-(v), it can then be shown (see Sirlanci et. al. [18,19]) that $\mathcal{A}^N \in G(M, \lambda_0)$ on $\mathcal{U}^N$ (i.e. that $|e^{\mathcal{A}^N t}|_{\mathcal{H}^N} \leq M e^{\lambda_0 t}, t \geq 0$), where the constants M and $\lambda_0$ are independent of $N$. Suppose further that for some $\lambda \geq \lambda_0$, $|\mathcal{J}^N R_\lambda(\mathcal{A})x - R_\lambda(\mathcal{A}^N)\mathcal{J}^N x|_{\mathcal{H}^N} \to 0$ as $N \to \infty$, for every $x \in \bar{\mathcal{H}}$, where $R_\lambda(\mathcal{A})$ and $R_\lambda(\mathcal{A}^N)$ denote respectively the resolvent operators $\mathcal{A}$ and $\mathcal{A}^N$ at $\lambda$. It then follows, from a version ([1,18,19]) of the Trotter Kato theorem (see, for example, [11,15]) that allows for the state spaces to depend on the parameters, $(\vec{a}^N, \vec{b}^N, \theta^N)$, that
$$|\mathcal{J}^N \mathcal{T}(t)x - e^{\mathcal{A}^N t}\mathcal{J}^N x|_{\mathcal{H}^N} \to 0 \text{ as } N \to \infty$$
for $x \in \bar{\mathcal{H}}$, uniformly in $t$ in compact intervals of $[0, \infty)$.

Define the operators $\mathcal{B}^N: \mathbb{R} \to \mathcal{U}^N$ and $\mathcal{C}^N: \mathcal{U}^N \to \mathbb{R}$ by
$$\langle \mathcal{B}^N u, v^N\rangle = \int_{Q^N} \langle b(q), v^N(q)\rangle_{V^*,V} f(q; \vec{a}^N, \vec{b}^N, \theta^N) dq u,$$
and
$$\mathcal{C}^N v^N = \int_{Q^N} \langle c(q), v^N(q)\rangle_{V^*,V} f(q; \vec{a}^N, \vec{b}^N, \theta^N) dq,$$

where $v^N \in \mathcal{U}^N$, and $u \in \mathbb{R}$ and consider the sequence of finite dimensional optimization problems given by

($\mathcal{P}_N$) Given $\nu$ data sets $\{(\tilde{u}_i, \tilde{y}_i)\}_{i=1}^{\nu}$, determine $\vec{a}^{N*} = [a_i^{N*}]$, $\vec{b}^{N*} = [b_i^{N*}]$ and $\theta^{N*} \in \Theta$, feasible, which minimize

$$J^N(\vec{a}, \vec{b}, \theta) = \sum_{i=1}^{\nu} \sum_{j=0}^{\mu_i} |y_j^N(\tilde{u}_i, \vec{a}, \vec{b}, \theta) - \tilde{y}_{i,j}|^2,$$

where $\{y_j^N(\tilde{u}_i, \vec{a}, \vec{b}, \theta)\}_{j=0}^{\mu_i}$ is given by

$$x_{i,j+1}^N = \hat{\mathcal{A}}^N x_{i,j}^N + \hat{\mathcal{B}}^N \tilde{u}_{i,j}, \quad y_{i,j}^N = \mathcal{C}^N x_{i,j}^N, \quad (4.1)$$
$$j = 0, 1, 2, \ldots \mu_i,$$

with $x_{i,0}^N = 0 \in \mathcal{U}^N$, where $\hat{\mathcal{A}}^N = e^{\mathcal{A}^N \tau} \in \mathcal{L}(\mathcal{U}^N, \mathcal{U}^N)$, $\hat{\mathcal{B}}^N = \int_0^{\tau} e^{\mathcal{A}^N s} \mathcal{B}^N ds \in \mathcal{L}(\mathbb{R}, \mathcal{U}^N)$, and $\hat{\mathcal{C}}^N = \mathcal{C}^N \in \mathcal{L}(\mathcal{U}^N, \mathbb{R})$. ∎

We can then prove the following convergence theorem (see [18,19]).

**Theorem 4.1** Let assumptions (i)-(v) hold. Then for each $N = 1, 2, \ldots$, the problems ($\mathcal{P}_N$) given above admit a solution $\hat{\rho}^N = (\vec{a}^{N*}, \vec{b}^{N*}, \theta^{N*})$. Suppose further that

vi. $|\mathcal{P}^N v^N - v^N|_{\mathcal{V}^N} \to 0$ as $N \to \infty$, $v^N \in \mathcal{V}^N$.

Then there exists a subsequence, $\{\hat{\rho}^{N_j}\}_{j=1}^{\infty}$, with $\hat{\rho}^{N_j} = (\vec{a}^{N_j*}, \vec{b}^{N_j*}, \theta^{N_j*}) \to (\vec{a}^*, \vec{b}^*, \theta^*) = \hat{\rho}$ as $j \to \infty$, and $\hat{\rho} = (\vec{a}^*, \vec{b}^*, \theta^*)$ a solution to problem ($\mathcal{P}$).

The optimization problems ($\mathcal{P}_N$) are solved numerically, typically via an iterative gradient-based scheme. Once a basis for the space $\mathcal{U}^N$ is chosen, the operators in (4.1) can be represented as matrices and the value of the cost functional $J$ and its gradient can then be computed. If for $N = 1, 2, \ldots$, $\mathcal{U}^N = \text{span}\{\psi_j^N\}_{j=1}^{K^N} \subset \mathcal{U}^N$, then the matrix representation for $\mathcal{A}^N \in \mathcal{L}(\mathcal{U}^N, \mathcal{U}^N)$ is given by $[\mathcal{A}^N]_{ij} = -[\langle \psi_i^N, \psi_j^N \rangle_{\mathcal{H}^N}]^{-1} [a(\psi_i^N, \psi_j^N)]$, for $i, j = 1, 2, \ldots, K^N$. Matrix representations for the operators $\mathcal{B}^N$ and $\mathcal{C}^N$ are computed analogously. The matrix representations for $\mathcal{A}^N$, $\mathcal{B}^N$, and $\mathcal{C}^N$ can then be used in a straight forward manner to compute the matrix representations for the operators $\hat{\mathcal{A}}^N$, $\hat{\mathcal{B}}^N$, and $\hat{\mathcal{C}}^N$ appearing in (4.1).

We compute $\vec{\nabla} J^N$ using the adjoint [12]. For each $i = 1, 2, \ldots, \nu$, set $v_{i,j}^N = [2(\hat{\mathcal{C}}^N x_{i,j}^N - \tilde{y}_{i,j}), 0, \ldots, 0]^T \in \mathbb{R}^{K^N}$, $j = 0, 1, 2, \ldots, \mu_i$ and define the adjoint systems

$$z_{i,j-1}^N = [\hat{\mathcal{A}}^N]^T z_{i,j}^N + v_{i,j-1}^N. \quad (4.2)$$

Then $\vec{\nabla} J^N$ at $(\vec{a}, \vec{b}, \theta)$ can then be computed as

$$\vec{\nabla} J^N(\rho) = \sum_{i=1}^{\nu} \sum_{j=1}^{\mu_i} [z_{i,j}^N]^T \left\{ \frac{\partial \hat{\mathcal{A}}^N}{\partial \rho} x_{i,j-1}^N \right.$$
$$-(\mathcal{A}^N)^{-1} \left\{ \frac{\partial \mathcal{A}^N}{\partial \rho} (\mathcal{A}^N)^{-1} (\hat{\mathcal{A}}^N - I) \frac{\partial \hat{\mathcal{B}}^N}{\partial \rho} \tilde{u}_{i,j-1} \right.$$
$$\left. - \frac{\partial \hat{\mathcal{A}}^N}{\partial \rho} \hat{\mathcal{B}}^N \tilde{u}_{i,j-1} - (\hat{\mathcal{A}}^N - I) \frac{\partial \hat{\mathcal{B}}^N}{\partial \rho} \tilde{u}_{i,j-1} \right\} \right\} \quad (4.3)$$
$$+ \sum_{i=1}^{\nu} \sum_{j=0}^{\mu_i} (y_j^N - \tilde{y}_{i,j})^T \frac{\partial \hat{\mathcal{C}}^N}{\partial \rho} x_{i,j}^N,$$

where $\rho = (\vec{a}, \vec{b}, \theta)$.

The tensor $\frac{\partial \hat{\mathcal{A}}^N}{\partial q}$ can be computed at the same time as the matrix $\hat{\mathcal{A}}^N$ is computed using the sensitivity equations. For $t \geq 0$ set $\Phi^N(t) = e^{\mathcal{A}^N t}$. Then $\Phi^N$ is the unique principal fundamental matrix solution to the initial value problem

$$\dot{\Phi}^N = \mathcal{A}^N \Phi^N, \Phi^N(0) = I. \quad (4.4)$$

Setting $\Psi^N = \partial \Phi^N / \partial \rho$, differentiating (4.4) with respect to $\rho$, interchanging the order of differentiation, and using the product rule, we obtain

$$\dot{\Psi}^N = \mathcal{A}^N \Psi^n + (\partial \mathcal{A}^N / \partial \rho) \Phi^N, \quad \Psi^N(0) = 0, \quad (4.5)$$

Combining the two initial value problems given in (4.4) and (4.5), and then solving we obtain

$$\begin{bmatrix} \frac{\partial \hat{\mathcal{A}}^N}{\partial \rho} \\ \hat{\mathcal{A}}^N \end{bmatrix} = \begin{bmatrix} \Psi^N(\tau) \\ \Phi^N(\tau) \end{bmatrix} = \exp\left( \begin{bmatrix} \mathcal{A}^N & (\partial \mathcal{A}^N / \partial \rho) \\ 0 & \mathcal{A}^N \end{bmatrix} \tau \right) \begin{bmatrix} 0 \\ I \end{bmatrix}. \quad (4.6)$$

## 5. Consistency of the Estimator

In the context of the alcohol biosensor problem of interest to us here, the estimator, $\hat{\rho}$ defined in (1.6) is given by, $\hat{\rho} = (\vec{a}^*, \vec{b}^*, \theta^*)$, where $\vec{a}^* = [a_i^*]$, $\vec{b}^* = [b_i^*]$ and $\theta^* \in \Theta$ are a solution to problem ($\mathcal{P}$) in Section 4.1. Under the following assumptions, using Theorem 4.2 in [4], it is possible to establish a consistency result for the estimator $\hat{\rho}$.

(a) The measurement noise $\{\varepsilon_{i,j}\}$ in (1.5) is i.i.d. with respect to a probability space $\{\Omega, \Sigma, P\}$ with $E_P[\varepsilon_{i,j}] = 0$ and $Var[\varepsilon_{i,j}] = \sigma^2$.
(b) The feasible set of parameters $\mathcal{R}$ is compact (i.e. closed and bounded since it is finite dimensional) and has nonempty interior.
(c) For $= 1, 2, \ldots, \nu$, $\mu_i = \mu$, and $\mu \tau = T$, for some positive integer $\mu$ and some $T > 0$, where $\tau$ is the sampling time defined in Section 3.
(d) Let $V = \{V_{i,j}\}$ be as is given in (1.5) for some $\rho_0 \in$ int $\mathcal{R}$ with $\tilde{y}_{i,j}(\tilde{u}_i, \rho_0)$ as given in the definition of problem ($\mathcal{P}$).
(e) For each $i = 1, 2, \ldots, \nu$, and $y(t; \tilde{u}_i, \rho) = y(t; \tilde{u}_i, \vec{a}, \vec{b}, \theta) = y(t)$ with $y(t)$ given by (3.5), $\rho_0$ is the unique minimizer of $J_{i,0}$ in $\mathcal{R}$, where

$$J_{i,0}(\rho) = \sigma^2 + \int_0^t (y(t; \tilde{u}_i, \rho_0) - y(t; \tilde{u}_i, \rho))^2 dt.$$

A straight forward application of Theorem 4.2 in [4] can then be used to establish the following lemma.

**Lemma 5.1** Assume that Assumptions (i)-(v) and (a) – (e) hold. Then there exists an $A \in \Sigma$ with $P(A) = 1$ such that for all $\omega \in A$ and $J_{i,\mu}(\rho; V)$ as given in (1.6) with $\alpha_{i,j} = 1$, we have

$$\frac{1}{\nu} \sum_{i=1}^{\nu} \left\{ \frac{1}{\mu} J_{i,\mu}(\rho; V(\omega)) - J_{i,0}(\rho) \right\} \to 0$$

as $\nu, \mu \to \infty, \tau \to 0$, with $\mu \tau = T$, uniformly in $\rho$, for $\rho \in \mathcal{R}$.

**Theorem 5.1** (Consistency of the estimator $\hat{\rho}$) Let $\hat{\rho} \in \mathcal{R}$ be as defined in (1.6) together with problem ($\mathcal{P}$). Then under the assumptions of Lemma 5.1, the estimator $\hat{\rho} = $ $(\vec{a}^*, \vec{b}^*, \theta^*)$ is consistent for $\rho_0$. That is $\hat{\rho} \xrightarrow{P} \rho_0$, as $\nu, \mu \to \infty, \tau \to 0$, with $\mu\tau = T$.

**Proof** The proof is quite similar to the proof of Theorem 4.3 in [4]. For $\nu = 1, 2, \ldots$, let $J_0(\rho) = \frac{1}{\nu}\sum_{i=1}^{\nu} J_{i,0}(\rho)$, let $A \in \Sigma$ be as in the statement of Lemma 5.1, let $\omega \in A$ be fixed, and let $\delta > 0$ be arbitrary. Then Lemma 5.1 implies that there exists an open local neighborhood of $\rho_0$ of radius $\delta$, $N_\delta(\rho_0)$, such that by Assumption (e) there exists an $\varepsilon > 0$ for which $J_0(\rho) - J_0(\rho_0) > \varepsilon$, for all $\rho \in \mathcal{R}_\delta$, where $\mathcal{R}_\delta$ is the compact (i.e. closed and bounded) set given by $\mathcal{R} \cap N_\delta(\rho_0)^c$. Now once again by Lemma 5.1, there exist $\nu_0, \mu_0, \tau_0$, such that for all $\nu > \nu_0, \mu > \mu_0, \tau < \tau_0$ with $\mu\tau = T$, $\left|\frac{1}{\nu\mu}J(\rho, V(\omega)) - J_0(\rho)\right| < \varepsilon/4$, $\rho \in \mathcal{R}$, where $J$ is as given in (1.6) or problem ($\mathcal{P}$). Then with $\nu > \nu_0, \mu > \mu_0, \tau < \tau_0, \mu\tau = T$ and $\rho \in \mathcal{R}_\delta$,

$$\frac{1}{\nu\mu}\big(J(\rho, V(\omega)) - J(\rho_0, V(\omega))\big)$$
$$= \frac{1}{\nu\mu}J(\rho, V(\omega)) - J_0(\rho) + J_0(\rho) - J_0(\rho_0) + J_0(\rho_0)$$
$$- \frac{1}{\nu\mu}J(\rho_0, V(\omega)) \geq -\frac{\varepsilon}{4} + \varepsilon - \frac{\varepsilon}{4} > 0.$$

But $J(\hat{\rho}, V(\omega)) \leq J(\rho_0, V(\omega))$. It follows that $\hat{\rho} \in N_\delta(\rho_0)$ if $\nu > \nu_0, \mu > \mu_0, \tau < \tau_0$, with $\mu\tau = T$. Since $\delta > 0$ was arbitrary and $P(A) = 1$, we have that $\hat{\rho} \xrightarrow{P} \rho_0$, as $\nu, \mu \to \infty, \tau \to 0$, with $\mu\tau = T$.

## 6. Numerical Results

The approximating finite dimensional subspaces $\mathcal{U}^N$ were constructed as follows based on the discretization of $(\eta, q_1, q_2)$-space. For each $n = 1, 2, \ldots$, let $\{\varphi_j^n\}_{j=0}^n$ denote the set of standard linear B-splines on the interval $[0,1]$ defined with respect to the usual uniform mesh, $\{j/n\}_{j=0}^n$, and set $H^n = \text{span}\{\varphi_j^n\}_{j=0}^n \subset V$ (note the $\varphi_j^n$ are the usual "pup tent" or "chapeau" functions of height one and support of width $2/n$, $\left[\frac{j-1}{n}, \frac{j+1}{n}\right] \cap [0,1]$). If $P^n : H \to H^n$ denotes the orthogonal projection of $H = L_2(0,1)$ onto $H^n$, it is well known (see for example, (Schultz, [21]) that $\lim_{n\to\infty} P^n\varphi = \varphi$ in $H$ for $\varphi \in H$ and in $V$ for $\varphi \in V$. Then for $i = 1, 2$, and each $m_i = 1, 2, \ldots$ let $\{\chi_{i,j}^{m_i}\}_{j=1}^{m_i}$ denote the set of standard $0^{\text{th}}$ order B-splines (i.e. piecewise constant functions) on the interval $[a_i, b_i]$ defined with respect to the usual uniform mesh, $\{a_i + (b_i - a_i)j/m_i\}_{j=0}^{m_i}$. If $P_i^{m_i}$ denotes the orthogonal projection of $L_2(a_i, b_i)$, it is not difficult to show [1,3] that $\lim_{m_i\to\infty} P_i^{m_i}\zeta = \zeta$ in $L_2(a_i, b_i)$ for every $\zeta \in L_2(a_i, b_i)$. Then let $N$ denote the triple $(n, m_1, m_2)$, and $L$ the multi-index $L = (j, j_1, j_2)$, where $j \in \{0, 1, 2, \ldots, n\}, j_1 \in \{1, 2, \ldots, m_1\}$, and $j_2 \in \{1, 2, \ldots, m_2\}$. We then use tensor products to define $\{\psi_L^N\}_L^N$ as $\psi_L^N = \varphi_j^n \chi_{1,j_1}^{m_1} \chi_{2,j_2}^{m_2}$ and set $\mathcal{U}^N = \text{span}\{\psi_L^N\}_L^N$. It is then not difficult to argue that Assumption (vi) holds.

In both of the examples to follow we fit a truncated bivariate normal. Let $\vec{q} = (q_1, q_2), \vec{a} = (a_1, a_2)$, and $\vec{b} = (b_1, b_2)$, and let $\phi(\cdot\,; \vec{\mu}, \Sigma)$ denote the joint density for the bivariate normal with mean $\vec{\mu}$ and covariance matrix $\Sigma$:

$$\phi(\vec{q}; \vec{\mu}, \Sigma) = \frac{1}{2\pi|\det\Sigma|^{1/2}}\exp\left(\frac{1}{2}(\vec{q}-\vec{\mu})\Sigma^{-1}(\vec{q}-\vec{\mu})^T\right).$$

Let $\Phi(\cdot\,; \vec{\mu}, \Sigma)$ denote the corresponding cumulative distribution function. We then set

$$f(\vec{q}; \vec{a}, \vec{b}, \vec{\mu}, \Sigma) = \frac{\phi(\vec{q}; \vec{\mu}, \Sigma)\chi_{[a_1,b_1]\times[a_2,b_2]}(\vec{q})}{\Phi(\vec{b}; \vec{\mu}, \Sigma) - \Phi((a_1,b_2); \vec{\mu}, \Sigma) - \Phi((b_1,a_2); \vec{\mu}, \Sigma) + \Phi(\vec{a}; \vec{\mu}, \Sigma)}.$$

In order to guarantee that we only search over positive definite symmetric matrices $\Sigma$, we parameterize $\Sigma$ as $\Sigma = LL^T$, where

$$L = \begin{bmatrix} L_{11} & 0 \\ L_{21} & L_{22} \end{bmatrix}.$$

It then follows that

$$\Sigma = LL^T = \begin{bmatrix} L_{11}^2 & L_{11}L_{21} \\ L_{11}L_{21} & L_{21}^2 + L_{22}^2 \end{bmatrix} > 0,$$

so long as $L_{11}$ and $L_{22}$ are both nonzero. Thus the optimization is over a feasible (e.g., since our model is diffusion based, we would want $a_1, a_2 > 0$ or $\bar{\alpha} > 0$) subset of $\mathbb{R}^9$ with $\rho = (a_1, b_1, a_2, b_2, \mu_1, \mu_2, L_{11}, L_{21}, L_{22})$.

All computations were carried out in Matlab on either MAC or PC laptops or desktops. For higher dimensional problems with high resolution discretization, faster platforms such as a cluster may be required. The optimization problems ($\mathcal{P}_N$) were solved iteratively using the Matlab Optimization Toolbox routine FMINCON. We computed the requisite gradients using the adjoint as shown in (4.2) – (4.6) above and we also let FMINCON compute them using finite differences. Both yielded the same results. The finite difference calculations were faster, but the adjoint would likely be preferable for problems involving a higher dimensional parameter space such as would be encountered in the non-parametric case. This is because the required number of integrations of the state equation when using the adjoint method does not increase with the number of parameters to be estimated as it does with a finite difference scheme for computing the gradient of the cost functional. Initial estimates for the parameters were obtained by first fitting each dataset deterministically via nonlinear least squares to obtain estimates for $q_1$ and $q_2$. Then sample means, standard deviations, and covariances were used to compute initial estimates for $a_1, b_1, a_2, b_2, \mu_1, \mu_2$. The two random variables $q_1$ and $q_2$ were initially assumed to be independent, each with standard deviation one sixth of the length of the corresponding boundary of the initial domain. Some care must be exercised in choosing these initial guesses. Because of the nature of the approximation scheme we are using, if in any iteration the pdf becomes too flat, our Galerkin scheme's mass and stiffness matrices can become singular or close to singular.

A WrisTAS™ 7 alcohol biosensor (Fig. 1.1) was worn for 18 days by one of the co-authors (S.E.L.) and was set to measure the local ethanol vapor concentration over the skin surface at 5-minute intervals. In addition, she contemporaneously collected breath measurements.

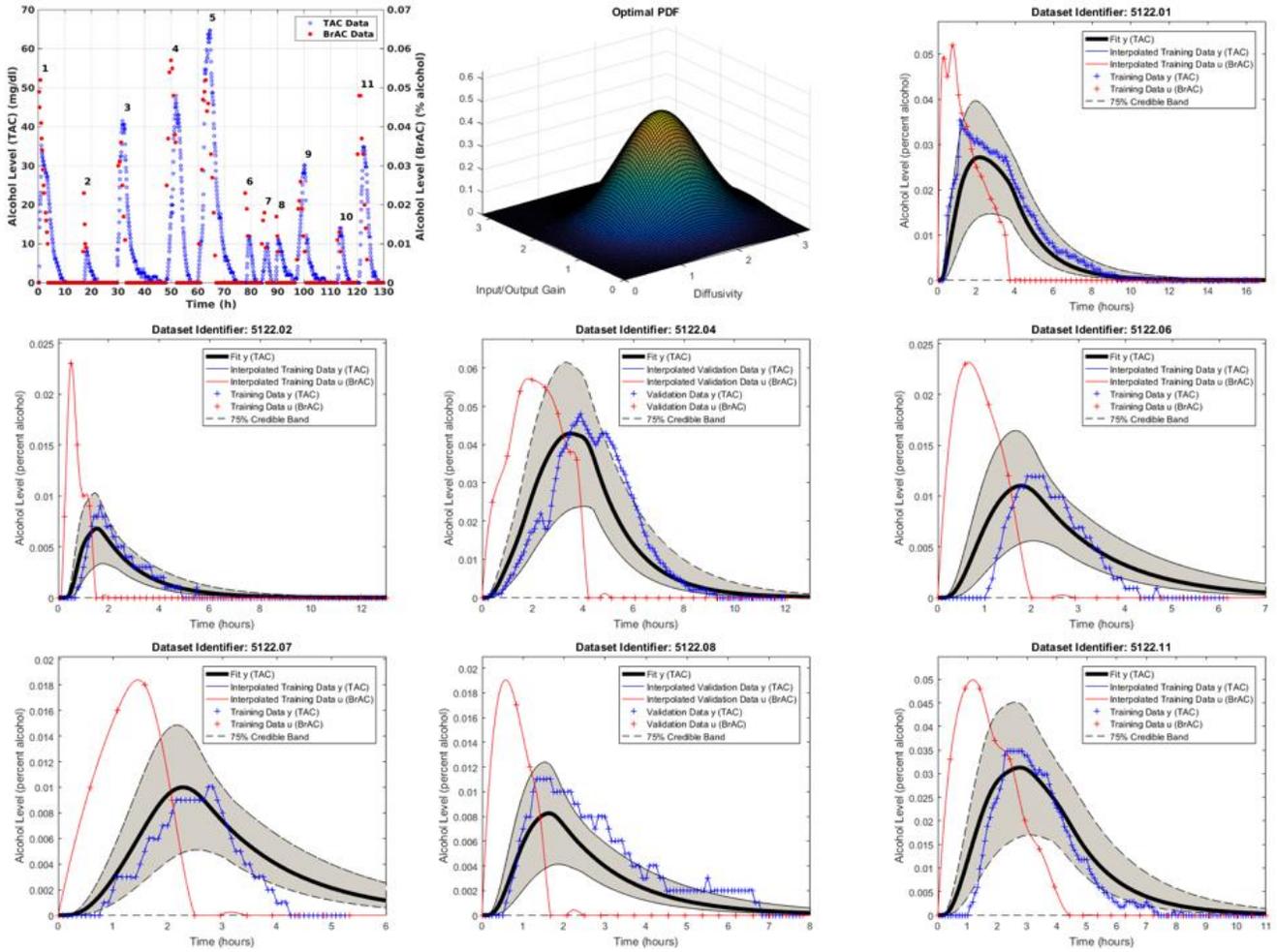

Figure 6.1. Row-wise from the top left: Panel (a) BrAC and TAC data for the 11 drinking episodes, Panel (b) the optimal pdf, Panels (c), (d), (f), (g), and (i) are the plots of training BrAC and corresponding TAC, resulting fit population model estimated TAC, and 75% credible band; Panels (e) and (h) are the cross validation results.

The first drinking episode was conducted in the laboratory with BrAC measured and recorded every 15 minutes from the start of the drinking session until BrAC returned to 0.000. She then wore the TAC device in the field and consumed alcohol ad libitum for the following 17 days. For each drinking episode, BrAC readings were taken every 30 minutes until the BrAC returned to 0.000. Figure 6.1.a shows the entire 18 day TAC signal along with the contemporaneous BrAC measurements. The 11 individual drinking episodes are marked. The TAC measurements provided by the sensor are in units of milligrams per deciliter (mg/dl), and the BrAC measurements are in units of percent alcohol.

We fit the population consisting of all eleven drinking episodes, but we also visually stratified the population into two groups, one containing the seven episodes in which the peak BrAC was higher than the (bench calibrated) peak TAC (episodes 1,2,4,6,7,8, and 11), and a second containing the remaining four drinking episodes in which the reverse was true (episodes 3, 5, 9, and 10). Our results for the first stratified group are shown in Figure 6.1.b - 6.1.i. In Figures 6.1.c,d,f,g, and i we plotted the training BrAC and TAC data for each of the episodes 1, 2, 6, 7, and 11 along with resulting fit population model estimated TAC and the 75% credible band. In Figures 6.1.e and 6.1.h we plotted the results of a cross validation on episodes 4 and 8. In Figure 6.1.b we plotted the optimal fit population truncated bivariate normal pdf. The converged values for the parameters were $a_1^* = 0.0000$, $b_1^* = 1.4850$, $a_2^* = 0.0000$, $b_2^* = 2.0363$, $\mu_1^* = 0.6318$, $\mu_2^* = 1.0295$, and $\Sigma^* = \begin{bmatrix} 0.0259 & 0.0077 \\ 0.0077 & 0.1232 \end{bmatrix}$.

The credible bands were computed directly from samples, $\vec{q} = (q_1, q_2)$ of $\vec{\tilde{q}} = (\tilde{q}_1, \tilde{q}_2)$, obtained from the optimal distribution for $\vec{\tilde{q}}$ using importance sampling and the state as $\langle c(\vec{q}), x_{i,j+1}^N \rangle_{V^*,V} = x_{i,j+1}^N(0, \vec{q})$, where $x_{i,j+1}^N = x_{i,j+1}^N(\eta, \vec{q})$ is given by (4.1). We note that, strictly speaking, this is not valid since our theory yields only that $x_{i,j+1}^N \in \mathcal{V} = L_\pi^2(Q; V)$, and thus that pointwise evaluation in $\vec{q}$ is undefined. However, the results seem quite reasonable and are extremely useful and consequently we have included them. Our results for the full un-stratified data set and for the second stratified group were similar. We also applied our scheme to a population consisting of multiple subjects each with a single drinking episode with the results being quite similar to those presented above [20].

## 7. Discussion and Concluding Remarks

We note that our theory and general approach also applies to hyperbolic systems with either H-semicoercive or V-coercive damping [2] such as the telegraph equation which can be used to model diffusion with finite speed of propagation. We are investigating elimination of the requirement that the measures $\pi$ be defined in terms of a density. We believe that it is possible to directly apply the Prohorov metric based framework developed in [4] by using a different version of the Trotter Kato-like semigroup approximation result in Section 4. We also believe that results for the estimation of functional parameters in parabolic systems could be used to estimate the pdfs non-parametrically. Indeed, in the system (3.2), (3.3), the pdf effectively plays the role of a non-constant coefficient in an abstract parabolic system. Thus, we should be able to estimate both the support and the shape of the density by parameterizing the pdf as a linear combination of basis elements (e.g. splines, orthogonal polynomials, etc.) and then estimating the coefficients in the expansion. We are also looking at polynomial chaos expansions for $\vec{q}$ and then estimating the coefficients.

Of primary interest to us is the estimation of the input to the system $u$, or the BAC/BrAC, from the output $y$, or TAC. Once the distribution of the random parameters has been estimated, this takes the form of a deconvolution problem. In [20] we use the results presented here together with the framework in [9] and [17] to do just that. We obtain an estimate of the input along with error bars or credible bands. Finally, we are looking at using the approach in [9] and [17] to control random parabolic systems, in particular, the computation of the feedback solution to the LQR and LQG problems for random distributed parameter systems.

**Acknowledgements.** This research was supported by a grants from the Alcoholic Beverage Medical Research Foundation and the National Institute of Alcohol Abuse and Alcoholism (NIAAA) (R21AA17711, R01AA026368-01 S.E.L. and I.G.R.), and (R01AA025969, C.E.F.).

**References**

[1] H.T. Banks, J.A. Burns and E.M. Cliff, Parameter estimation and identification for systems with delays, SIAM J. Contr. and Opt., Vol. 19, No. 6, November 1981, pp. 791-828.

[2] H. T. Banks. and K. Ito, A Unified framework for approximation in inverse problems for distributed parameter systems, Control Theory Advanced Technology Vol. 4, No. 1, 1988, pp. 73-90.

[3] H. T. Banks, and K. Kunisch, Estimation techniques for distributed parameter systems, Boston, Birkhauser, 1989.

[4] H.T. Banks and C. Thompson, Least squares estimation of probability measures in the Prohorov metric framework, CRSC-TR12-21, N. C. State University, 2012.

[5] Z. Dai, I.G. Rosen, C.Wang, N.Barnett, and S.E. Luczak, Using drinking data and pharmacokinetic modeling to calibrate transport model and blind deconvolution based data analysis software for transdermal alcohol biosensors, Math. Biosci. and Eng., 13(5), 2016, pp. 911–934.

[6] M. Davidian, and D. Giltinan, Nonlinear Models for Repeated Measurement Data,NewYork: Chapman and Hall, 1995.

[7] E. Demidenko, Mixred Models, Theory and Applications, John Wiley and Sons, Hoboken, 2004.

[8] C. E. Fairbairn, K. Bresin, D. Kang, I. G. Rosen, T. Ariss, N. P. Barnett, S. E. Luczak, and N. S. Eckland, A multimodal investigation of contextual effects on alcohol's emotional rewards, 2017, submitted, in review.

[9] C. J. Gittelson, R. Andreev, and C. Schwab, Optimally adaptive Galerkin methods for random partial differential equations, Journal of Computational and Applied Mathematics, 263, 2014, pp.189–201.

[10] E. Grenier, V. Louvet, and P. Vigneaux, Parameter estimation in non-linear mixed effects models with SAEM algorithm: extension from ODE to PDE. ESAIM: Math. Modelling and Num. Anal., 48(5), 2014, pp. 1303-1329.

[11] T. Kato, Perturbation Theory for Linear Operators, Springer, 1976.

[12] A. Levi and I. G. Rosen, A novel formulation of the adjoint method in the optimal design of quantum electronic devices. SIAM J. Ctrl and Opt, 48, 2010, pp. 3191–3223.

[13] J. L. Lions, Optimal Control of Systems Governed by Partial Differential Equations, New York: Springer, 1971.

[14] D.R. Mould and R.N. Upton, Basic concepts in population modeling, simulation, and model-based drug development, CPT: Pharmacomet. & Systems Pharm. 2012, 1, e6.

[15] A. Pazy, Semigroups of Linear Operator and Applications to PDEs, New York Springer, 1983.

[16] I.G. Rosen, S.E. Luczak, and J. Weiss, Blind deconvolution for distributed parameter systems with unbounded input and output and determining blood alcohol concentration from transdermal biosensor data. Applied Mathematics and Computation, 231, 2014, pp. 357–376.

[17] C. Schwab and C. J. Gittelson, Sparse tensor discretization of high-dimensional parametric and stochastic PDEs, in Acta Num, 20, Cambridge U. Press, 2011, pp. 291-467.

[18] M. Sirlanci, S. E. Luczak and I. G. Rosen, Approximation and convergence in the estimation of random parameters in linear holomorphic semigroups generated by regularly dissipative operators, Proceedings of the 2017 American Control Conference, May 2017.

[19] M. Sirlanci and I. G. Rosen, Estimation of the distribution of random parameters in discrete time abstract parabolic systems with unbounded input and output: approximation and convergence, J. Math Anal and App., Submitted, 2017. Retrieved from https://arxiv.org/pdf/1804.04904.pdf

[20] M. Sirlanci, S. E. Luczak, C. E. Fairbairn, K. Bresin, D. Kang, and I. G. Rosen, Deconvolving the input to random abstract parabolic systems; a population model-based approach to estimating blood/breath alcohol concentration from transdermal alcohol biosensor data, submitted, 2018. Retrieved from https://arxiv.org/pdf/1807.05088.pdf

[21] M.H. Schultz, Spline Analysis, Prentice Hall, Englewood Cliffs, N.J., 1973.

[22] H. Tanabe, Equations of Evolution (Vol. 6). Pitman Publishing, 1979.